\documentclass[a4paper,12pt]{article}

\usepackage{amssymb}
\usepackage{tipa
}

\usepackage[T2A]{fontenc}
\usepackage[cp1251]{inputenc}
\usepackage{amsmath}

\usepackage[dvips]{graphicx}
\usepackage[
english,
russian%,%for the right language of the date
]{babel}

\usepackage{wasysym}
\usepackage{bbm}

\def \GK{\operatorname{GK}}

\def\Ht{\operatorname{Ht}}
\def\PI{{\rm PI}}

\newtheorem{theorem}{Теорема}[section]

\newtheorem{lemma}[theorem]{Лемма}
\newtheorem{proposition}[theorem]{Предложение}
\newtheorem{notation}[theorem]{Обозначение}
\newtheorem{remark}[theorem]{Замечание}
\newtheorem{construction}[theorem]{Конструкция}

\newtheorem{definitionhead}[theorem]{Определение}
\newenvironment{definition}{\begin{definitionhead}%
\sl}{\end{definitionhead}}

\begin{document}

%\title{Субэкспоненциальная оценка  в теореме Ширшова о высоте}

\title{
\ \hbox to \textwidth{\normalsize      % Dirty trick to shift
УДК  512.552.4+512.57+519.1\hfill}\\[1ex]   % УДК to the left
Estimations of the particular periodicity in case of the extremal periods in Shirshov's Height theorem}

\author{Mikhail Kharitonov}

\maketitle

\begin{abstract}

Gelfand-Kirillov dimension of l-generated general matrixes is $(l-1)n^2+1.$ The minimal degree of the identity of this algebra is $2n$ as a corollary of Amitzur-Levitsky theorem. That is why essential height of $A$ --- $l$-generated PI-algebra of degree n --- over every set of words can be bigger than $(l-1)n^2/4 + 1.$ We prove in this paper, that if $A$ has finite $\GK$-dimension, then the number of lexicographically comparable subwords

1) with period 2 in every monoid of $A$ is not bigger than ${(2l-1)(n-1)(n-2)\over 2};$

2) with period 3 in every monoid of $A$ is not bigger than $(2l-1)(n-1)(n-2);$

3) with period $(n-1)$ in every monoid of $A$ is not bigger than $(l-2)(n-1).$

The cases of the words of length 2 and 3 is synthesized to the proof of Shirshov's Height theorem.

{\bf Key words:} Essential height, Shirshov's Height theorem, combinatorics of words, $n$-divisibility, Dilworth's theorem.

\end{abstract}

\section{Введение.}

В 1958 году А.~И.~Ширшов
доказал свою знаменитую теорему о высоте (\cite{Shirshov1}, \cite{Shirshov2}).
\begin{definition} \label{DfHeight}
Назовём \PI-алгебру $A$ алгеброй {\em ограниченной высоты
$\Ht_Y(A)$ над множеством слов $Y = \{ u_1, u_2,\ldots\}$}, если
$\Ht_Y(A)$ --  такое минимальное число, что любое слово $x$ из $A$ можно
представить в виде
$$x = \sum_i \alpha_i u_{j_{(i,1)}}^{k_{(i,1)}}
u_{j_{(i, 2)}}^{k_{(i,2)}}\cdots
u_{j_{(i,r_i)}}^{k_{(i,r_i)}},
$$
причем $\{r_i\}$ ограничены числом $\Ht_Y(A)$ в совокупности. Множество
$Y$ называется {\em базисом Ширшова} для $A$.
\end{definition}

\begin{definition}    \label{Dfndivided}
Слово $W$ называется {\em $n$-разбиваемым}, если его можно представить в
виде $W=W_0W_1\cdots W_n$, где подслова $W_1,\dots,W_n$ идут в
порядке лек\-си\-ко\-гра\-фи\-чес\-ко\-го убывания.
\end{definition}

\medskip
 {\bf Теорема Ширшова о высоте.} (\cite{Shirshov1}, \cite{Shirshov2})
 {\it Множество всех не $n$-разбиваемых слов в конечно порождённой алгебре с допустимым по\-ли\-но\-ми\-аль\-ным тождеством имеет ограниченную
 высоту $H$ над множеством слов степени не выше $n-1$.
 }
\medskip

Интерес к проблеме не ослабевает, возникают смежные вопросы и попытки их решения: см.,
например, \cite{BBL97,BK,Bel07,Kem09,Ch01,BP07,Lot83,02}.
 В работе \cite{BK} было получено, что высота в условиях теоремы Ширшова
 меньше $$\Phi(l,n) = 2^{87} l\cdot n^{12\log_3 n + 48}.$$

В данной работе исследуются возможности получения точной оценки.
Для этого ведётся исследование в двух направлениях:
\begin{enumerate}
    \item Поиск наименьшей верхней оценки.
    \item Поиск не $n$-разбиваемого слова как можно большей высоты.
\end{enumerate}

Для улучшения оценок в теореме Ширшова о высоте, полученных ранее в работе
\cite{BK}, необходимо оценить выборочную высоту над множествами
нециклических слов определённой длины. Мы рассматриваем случай,
когда означенная длина равна $2$, $3$ и $(n-1)$. Это имеет и
самостоятельную ценность, так как к нему можно свести
доказательство теоремы Ширшова с помощью кодировки.
\medskip

{\bf Основные результаты} работы состоят в следующем:
\begin{theorem}     \label{verh}
Малая выборочная высота множества не сильно $n$-разбиваемых слов
над $l$-буквенным алфавитом относительно множества нециклических
слов длины $2$ не больше $\beth(2,l,n)$, где
$$\beth(2,l,n) = {(2l-1)(n-1)(n-2)\over 2}.$$
\end{theorem}

\begin{theorem} \label{niz}
Малая выборочная высота множества не сильно $n$-разбиваемых слов
над $l$-буквенным алфавитом относительно множества нециклических
слов длины $2$ при фиксированном $n$ больше, чем $\Psi(n,l)$, где
$$\Psi(n,l) =  {n^2 l\over 2} (1 - o(l)).$$
Более точно, $\Psi(n,l) =(l-2^{n-1})(n-2)(n-3)/2.$
\end{theorem}

Теорема \ref{verh} с помощью кодировки обобщается до доказательства экс\-по\-нен\-ци\-аль\-ной оценки на существенную высоту:

\begin{theorem} \label{ess:t}
Существенная высота $l$-порождённой $PI$-алгебры с допустимым по\-ли\-но\-ми\-аль\-ным тождеством степени $n$ над множеством слов длины $<n$ меньше, чем $\Upsilon (n, l),$ где
$$\Upsilon (n, l) = 8(l+1)^n n^6l.$$
\end{theorem}

В работе \cite{BK} высота свободной $l$-порождённой алгебры $A$ с тождеством степени $n$ оценивается суммой $\Psi(n,4n,l)$ и существенной высоты алгебры $A$, где $\Psi(n,d,l)=2^{18} l (nd)^{3 \log_3 (nd)+13}d^2.$ Таким образом, из теоремы \ref{ess:t} вытекает следующая

\begin{theorem}
Существенная высота $l$-порождённой $PI$-алгебры с допустимым по\-ли\-но\-ми\-аль\-ным тождеством степени $n$ над множеством слов длины $<n$ меньше, чем $\Phi (n, l),$ где
$$\Phi (n, l) = 8(l+1)^{n+1} n^6.$$
\end{theorem}

Для исследования общего случая полезно оценить длину кусочно-периодических слов с разными длинами периодов.

 \begin{theorem}     \label{verh1}
Малая выборочная высота множества не сильно $n$-разбиваемых слов
над $l$-буквенным алфавитом относительно множества нециклических
слов длины $3$ не больше $\beth(3,l,n)$, где
$$\beth(3,l,n) = {(2l-1)(n-1)(n-2)}.$$
\end{theorem}

 \begin{theorem}     \label{verh2}
Малая выборочная высота множества не сильно $n$-разбиваемых слов
над $l$-буквенным алфавитом относительно множества нециклических
слов длины $(n-1)$ не больше $\beth(n-1,l,n)$, где
$$\beth(n-1,l,n) = (l-2)(n-1).$$
\end{theorem}

Малую и большую выборочные высоты связывает следующая

\begin{theorem} \label{co1}
Большая выборочная высота $l$-порождённой $\PI$-алгебры с
до\-пус\-ти\-мым полиномиальным тождеством степени $n$ над множеством
нециклических слов длины $k$ меньше $2(n - 1)\beth(k,l,n)$.
\end{theorem}

Данная оценка практически полученную А.~Я.~Беловым в \cite{Bel92}, но получена другими методами, в частности с использованием {\it графов Рози}. Если $W$ -- слово над алфавитом $\textsuperimposetilde{A}$, то {\it $k$-графом Рози} называется граф, вершины которого соответствуют различным подсловам слова $W$ длины $k$. Из вершины $w_1$ в вершину $w_2$ ведёт стрелка, если максимальный суффикс $w_1$ совпадает с максимальным префиксом $w_2$, то есть $w_1=a_1u, w_2=ua_2,$ где $a_1,a_2\in \textsuperimposetilde{A}.$

Таким образом, слову $W$ отвечает некоторая траектория в $k$-графе Рози. Фик\-си\-ру\-ем натуральное число $k$. Тогда теорему \ref{ess:t} можно переписать следующим образом:

\begin{theorem}
Пусть $\textsuperimposetilde{A}$ -- алфавит мощности $l$. Для любого натурального числа $n$ найдутся числа $N=N(n),$ $ d=d(n)$ такие, что на траектории любого не $n$-разбиваемого слова $W$ над алфавитом $\textsuperimposetilde{A}$ в $k$-графе Рози найдётся менее, чем $\Upsilon (n, l)$ циклов длины $<n$, по которым траектория проходит больше $d(n)$ раз.
\end{theorem}

{\bf О нижних оценках.} Сопоставим полученные результаты  с нижней
оценкой для высоты. Высота алгебры $A$ не меньше ее размерности
Гель\-фан\-да--Ки\-рил\-ло\-ва $\GK(A)$. Для алгебры
$l$-порождённых общих матриц порядка $n$ данная размерность равна
$(l-1)n^2+1$ (см. \cite{Procesi}, а также \cite{Bel04}). В то же
время, минимальная степень тождества этой алгебры равна $2n$ в
силу теоремы Ам\-и\-цу\-ра--Ле\-виц\-ко\-го. Имеет место следующее

\begin{proposition}
Высота $l$-порожденной $\PI$-алгебры степени $m$, а также
мно\-жест\-ва  не $m$-раз\-би\-ва\-ем\-ых слов над $l$-бук\-вен\-ным
алфавитом не менее, чем $(l-1)m^2/4+1$.
\end{proposition}

Эта оценка линейна по числу образующих $l$.

Нижние оценки на индекс нильпотентности были установлены
Е.~Н.~Кузьминым в работе \cite{Kuz75}. Он привел
пример $2$-порожденной алгебры с тождеством $x^n=0$, индекс
нильпотентности которой строго больше $(n^2+n-2)/2$. В то же время
для случая нулевой характеристики и счетного числа образующих
Ю.~П.~Размыслов (см., например, \cite{Razmyslov3}) получил верхнюю
оценку на индекс нильпотентности равную $n^2$.
\medskip

{\bf Благодарности.} Автор признателен В.~Н.~Латышеву,
А.~В.~Михалёву, участникам семинара ``Теория колец'' и участникам семинара МФТИ под руководством А.~М.~Рай\-го\-род\-ско\-го за
внимание к работе и лично А.~Я.~Белову за крайне полезные обсуждения.

\section{Основные понятия.}

\begin{definition}
Слово $u$ назовем {\em нециклическим}, если $u$ нельзя представить
в виде $v^k$, где $k>1$, $k\in \mathbb N.$ Слова $u$ и $v$ {\em несравнимы}, если
одно из них является началом другого. Слова $u$ и $v$ назовём {\em
сильно сравнимыми}, если любые их циклические сдвиги сравнимы.
\end{definition}

\begin{notation}
Пусть $a_1, a_2,\ldots,a_l$ -- буквы алфавита. {\em Значение}
буквы $a_i$ счи\-та\-ем равным $i$. $a_i$ лексикографически больше
$a_j,$ если $i>j$.
\end{notation}

Заметим, что понятие {\it сильной сравнимости} разбивает слова на
классы эк\-ви\-ва\-лент\-нос\-ти.

\begin{definition}
{\em Слово-цикл} -- некоторое слово $u$ со всеми его
сдвигами по циклу.
\end{definition}

\begin{definition}
Алгебра $A$ имеет {\em существенную высоту $h=H_{Ess}(A)$} над
ко\-неч\-ным множеством $Y$, называемым {\em $s$-базисом алгебры $A$},
если можно выбрать такое ко\-неч\-ное множество $D\subset A$, что $A$
линейно представима элементами вида $t_1\cdot\ldots\cdot t_l$, где
$l\leqslant 2h+1$, и $\forall i (t_i\!\in\! D \vee
t_i=y_i^{k_i};y_i\in Y)$, причем множество таких $i$, что
$t_i\not\in D$, содержит не более $h$ элементов. Аналогично
определяется {\em существенная высота} множества слов.
\end{definition}

\begin{definition}
а) Число $h$ называется {\em малой выборочной высотой} с границей $k$
слова $W$ над множеством слов $Z$, если $h$ -- такое максимальное
число, что у слова $W$ найдётся $h$ попарно непересекающихся
циклически несравнимых подслов вида $z^m,$ где $z\in Z, m>k$.

б) Число $h$ называется {\em большой выборочной высотой} с границей $k$
слова $W$ над множеством слов $Z$, если $h$ -- такое максимальное
число,что у слова $W$ найдётся $h$ попарно непересекающихся
подслов вида $z^m,$ где $z\in Z, m>k$, причём соседние подслова из
этой выборки несравнимы.

{\bf Здесь и далее:} $k = 2n$.

г) Слово $W$ называется {\em сильно $n$-разбиваемым}, если его
можно представить в виде $W=W_0W_1\cdots W_n$, где подслова
$W_1,\dots,W_n$ идут в порядке лексикографического убывания, и
каждое из слов $W_i, i=1, 2,\ldots, n$ начинается с некоторого
слова $z_i^k\in Z$ где все $z_i$ различны.

д) Множество слов $V$ имеет малую (большую) выборочную высоту $h$
над некоторым множеством слов $Z$, если $h$ является точной
верхней гранью малых (больших) вы\-бо\-роч\-ных высот над $Z$ его
элементов.
\end{definition}

Говоря неформально, любое длинное слово есть произведение
периодических час\-тей и ``прокладок'' ограниченной длины. {\em
Существенная высота} есть число таких пе\-ри\-о\-ди\-чес\-ких кусков, а {\em
выборочная высота} учитывает только куски определённого вида.

$s$-базис является базисом Ширшова тогда и только тогда, когда он
порождает $A$ как алгебру. Связь существенной высоты и размерности
Гельфанда-Кириллова рассмотрена в работе \cite{BBL97}.

{\bf $n$-разбиваемость и теорема Дилуорса.} Значение понятия {\it
$n$-раз\-би\-ва\-ем\-ос\-ти} выходит за рамки проблематики,
относящейся к проблемам бернсайдовского типа. Оно играет роль и
при изучении полилинейных слов, в оценке их количества, где {\it
полилинейным} называется слово, в которое каждая буква входит
не более одного раза. В.~Н.~Латышев применил теорему Дилуорса для
получения оценки числа не являющихся $m$-разбиваемыми полилинейных слов степени $n$ над
алфавитом\linebreak $\{a_1,\dots,a_n\}$. Эта
оценка:  ${(m - 1)}^{2n}$, и она близка к реальности. Напомним эту
теорему.

\medskip
{\bf Теорема Дилуорса}: {\it Пусть $n$ -- наибольшее количество
элементов ан\-ти\-це\-пи данного конечного частично упорядоченного
множества $M$. Тогда $M$ можно раз\-бить на $n$  попарно
непересекающихся цепей.}
\medskip

Впервые предложил применить эту идею к неполилинейному случаю
Г.~Р.~Челноков в 1996 году.

%Рассмотрим полилинейное слово $W$ из $n$ букв. Положим $a_i\succ
%a_j$, если $i>j$ и буква $a_i$ стоит в слове $W$ правее $a_j$.
%Условие не $k$-разбиваемости означает отсутствие антицепи из $n$
%элементов. По теореме Дилуорса тогда все позиции (и,
%соответственно, буквы $a_i$) разбиваются на $n-1$ цепь. Сопоставим
%каждой цепи свой цвет. Тогда раскраска позиций и раскраска букв
%однозначно определяет слово $W$. А число таких раскрасок не
%превосходит $(n-1)^k\times (n-1)^k=(n-1)^{2k}$.

\section{Доказательство оценок.}

Далее будем считать, что слова строятся над алфавитом $\textsuperimposetilde{A}$ из букв $\{a_1, a_2,\ldots,a_{l}\},$ над которыми введён лексикографический порядок, причём $a_i<a_j$, если $i<j$. Для следующих ниже доказательств будем отждествлять буквы $a_i$ с их индексами $i$ (то есть будем писать не слово $a_ia_j$, слово $ij$).

\subsection{Доказательство теоремы \ref{verh}% и \ref{verh1}
}

Пусть слово $W$ не сильно $n$-разбиваемо.  Рассмотрим некоторое
множество $\Omega$ попарно непересекающихся циклических несравнимых подслов
$W$ вида $z^m$, где $m>2n$, $z$~---~нециклическое двубуквенное
слово. Будем называть элементы этого множества {\it представителями,}
имея ввиду, что эти элементы являются представителями раз\-лич\-ных
классов эквивалентности по сильной сравнимости. Пусть набралось
$t$ таких пред\-ста\-ви\-те\-лей. Пронумеруем их всех в порядке положения
в слове $W$ (первое -- самое левое) числами от $1$ до $t$. В
каждом выбранном представителе в качестве подслов содержатся ровно
два различных двубуквенных слова.

Введём порядок на этих словах следующим образом: $u\prec v$, если
\begin{itemize}
    \item $u$ лексикографически меньше $v$,
    \item представитель, содержащий $u$
левее представителя, содержащего $v$.
\end{itemize}
Из не сильной
$n$-разбиваемости получаем, что максимальное возможное число
попарно несравнимых элементов равно $n-1$. По теореме Дилуорса
существует разбиение рассматриваемых двубуквенных слов на $(n-1)$
цепь. Раскрасим слова в $(n-1)$ цвет в соответствии с их
принадлежностью к цепям.

Введём соответствие между следующими четырьмя
объектами:
\begin{itemize}
    \item натуральными числами от $1$ до $t$,
    \item классами эквивалентности по сильной сравнимости,
    \item содержащимися в классах эквивалентности по сильной сравнимости цик\-ли\-чес\-ки\-ми словами
длины $2$,
    \item парами цветов, в которые
раскрашены сдвиги по циклу этого слово-цикла.
\end{itemize}

Буквы слово-цикла раскрасим в цвета, в которые раскрашены сдвиги по циклу, которые с них начинаются.

Рассмотрим граф $\Gamma$ с вершинами вида $(k,i)$, где $0<k<n,$ $0<i
\leqslant l$. Первая координата соответствует цвету, а вторая -- букве. Две вершины $(k_1,i_1), (k_2,i_2)$ со\-е\-ди\-ня\-ют\-ся
{\it ребром с весом  $j,$} если
\begin{itemize}
    \item в $j$-ом представителе содержится
слово-цикл из букв $i_1, i_2$,
    \item буквы $j$-го представителя раскрашены в
цвета $k_1, k_2$ соответственно.
\end{itemize}

В нашем графе будет не менее
$(n-1)$-ой компоненты связности, так как слова состоят из букв
разного цвета.

Посчитаем число рёбер между вершинами вида $(k_1,i_1)$ и вершинами
вида $(k_2,i_2)$, где $k_1, k_2$ -- фиксированы, $i_1, i_2$ --
произвольны. Рассмотрим два ребра $l_1$ и $l_2$ из рассматриваемого множества с
весами $j_1 < j_2$ с концами в некоторых вершинах $A=(k_1, i_{1_1}),
B=(k_2,i_{2_1})$ и $C=(k_1,i_{1_2}), D=(k_2,i_{2_2})$ соответственно.
Тогда по построению одновременно выполняются неравенства $ i_{1_1}\leq i_{1_2}, i_{2_1}\leq  i_{2_2}.$
При этом, так как рассматриваются представители классов
эквивалентности по сильной срав\-ни\-мос\-ти, то одно из неравенств
строгое. Значит, $i_{1_1}+i_{2_1}\leqslant i_{1_2}+i_{2_2}+1.$
Так как вторые координаты вершин ограничены числом $l$,
то вычисляемое число рёбер будет не более $(2l-1)$.

Так как первая координата вершин меньше $n$, то всего рёбер в
графе будет не более $(2l-1)(n-1)(n-2)/2$. Таким образом, теорема
\ref{verh} доказана.

\subsection{Доказательство теоремы \ref{verh1}}

%{\bf Докажем  теорему \ref{verh1}.}
Пусть слово $W$ не сильно
$n$-разбиваемо.  Рассмотрим некоторое множество попарно
непересекающихся циклических несравнимых подслов слова $W$ вида $z^m$,
где $m>2n$, $z$ -- нециклическое трёхбуквенное слово. Будем называть
элементы этого множества представителями, имея ввиду, что эти
элементы являются представителями различных классов
эквивалентности по сильной сравнимости. Пусть набралось $t$
таких пред\-ста\-ви\-те\-лей. Пронумеруем их всех в порядке положения в
слове $W$ (первое -- ближе всех к началу слова) числами от $1$ до $t$. В каждом
выбранном представителе в качестве подслов содержатся ровно три
различных трёхбуквенных слова.

Введём порядок на этих словах следующим образом:

$u\prec v$, если
\begin{itemize}
    \item $u$ лексикографически меньше $v$,
    \item представитель, содержащий $u$,
левее представителя, содержащего $v$.
\end{itemize}
 Из не сильной
$n$-разбиваемости получаем, что максимальное возможное число
попарно несравнимых элементов равно $n-1$. По теореме Дилуорса
существует разбиение рассматриваемых трёхбуквенных слов на $(n -
1)$ цепь. Раскрасим слова в $(n-1)$ цвет в соответствии с их
принадлежностью к цепям.

Рассмотрим теперь уже ориентированный граф $G$ с вершинами вида $(k,i)$, где
$0 < k < n,$ $0 < i \leqslant l$. Первая координата обозначает цвет, а вторая -- букву.
{\it Ребро с некоторым весом $j$}
выходит из $(k_1, i_1)$ в $(k_2, i_2)$, если для некоторого $i_3$
\begin{itemize}
    \item в $j$-ом представителе содержится слово-цикл $i_1 i_2 i_3$,
    \item буквы $i_1, i_2$ $j$-го представителя раскрашены в цвета $k_1, k_2$ соответственно.
\end{itemize}
Таким образом, граф $G$ состоит из ориентированных треугольников с
рёбрами одинакового веса. Однако, в отличие от графа $\Gamma$ из доказательства теоремы \ref{verh}, могут
появляться кратные рёбра, то есть рёбра с общими началом и концом, но разным весом. Для дальнейшего доказательства нам
потребуется

\begin{lemma}[Основная]
Пусть $A$, $B$ и $C$ -- вершины графа $G$, $A\to B\to C\to A$ --
ориентированный треугольник с рёбрами некоторого веса $j$, кроме того,
су\-щест\-ву\-ют другие рёбра $A\to B, B\to C, C\to A$ с весами $a, b,
c$ соответственно. Тогда среди $a, b, c$ есть число, большее $j$.
\end{lemma}

$\RHD$ От противного. Если два числа из набора $a, b, c$ равны друг другу,
то  $a = b = c = j$, так как в противном случае есть 2 треугольника $A\to B\to C\to A$, в каждом из которых веса всех трёх рёбер совпадают между собой. Тогда в $\Omega$ есть два не сильно сравнимых слова, что противоречит определению $\Omega.$ Без ограничения общности, что $a$ наибольшее из чисел
$a, b и с$. Рассмотрим треугольник из рёбер веса $a.$ Этот треугольник будет иметь общую с $\triangle ABC$ сторону $AB$ и некоторую третью вершину $C'.$ Если вторая координата точки $C'$ совпадает со второй координатой точки $C$ (то есть совпали соответствующие $C$ и $C'$ буквы алфавита), то $\triangle ABC$ и $\triangle ABC'$ соответствуют не сильно сравнимым словам из множества $\Omega.$ Снова получено противоречие с определением множества $\Omega.$ По предположению $a<j,$ а значит, из монотонности цвета $k_A$ (первой координаты вершины $A$) слово $i_A i_B i_{C'},$ составленное из вторых координат вершин $A, B, C'$ соответственно, лексикографически меньше слова $i_A i_B i_C.$
Значит, $i_{C'} < i_C.$ Тогда слово $i_B i_{C'}$ лексикографически меньше слова $i_B i_{C}.$ Из монотонности цвета $k_B$ получаем, что $b>a.$ Противоречие.
$\LHD$
\medskip

{\bf Завершение доказательства теоремы \ref{verh1} }

Рассмотрим теперь граф $G_1$, полученный из графа $G$ заменой
между каждыми двумя вершинами кратных рёбер на ребро с наименьшим
весом. Тогда в графе $G_1$ встретятся рёбра всех весов от 1 до
$t$.

Посчитаем число рёбер из вершин вида $(k_1, i_1)$ в вершины вида
$(k_2, i_2)$, где $k_1, k_2$ фиксированы, $i_1, i_2$
произвольны. Рассмотрим два ребра из рассматриваемого мно\-жест\-ва с
весами $j_1 < j_2$ с концами в некоторых вершинах $(k_1, i_{1_1}),
(k_2, i_{2_1})$ и $(k_1,i_{1_2}),(k_2,i_{2_2})$ соответственно.
Тогда по построению $ i_{1_1}\leq i_{1_2}, i_{2_1}\leq i_{2_2}$,
причём, так как рас\-смат\-ри\-ва\-ют\-ся представители классов
эквивалентности по сильной сравнимости, то одно из неравенств
строгое. Так как вторые координаты вершин ограничены числом $l$,
то вычисляемое число рёбер будет не более $2l-1$.

Так как первая координата вершин меньше $n$, то всего рёбер в
графе будет не более $(2l-1)(n-1)(n-2)$. Таким образом, теорема
\ref{verh1} доказана.

\subsection{Доказательство теоремы \ref{verh2}}

Пусть слово $W$ не $n$-разбиваемо. Как и прежде, рассмотрим некоторое множество попарно непересекающихся несравнимых подслов слова $W$ вида $z^m,$ где $m>2n,$ $z$ -- $(n-1)$-буквенное нециклическое слово.  Будем называть
элементы этого множества {\it представителями}, имея ввиду, что эти
элементы являются представителями раз\-лич\-ных классов
эквивалентности по сильной сравнимости. Пусть набралось $t$
таких предс\-та\-ви\-те\-лей. Пронумеруем их всех в порядке положения в
слове $W$ (первое -- ближе всех к началу слова) числами от $1$ до $t$. В каждом выбранном представителе в качестве подслов содержатся ровно $(n-1)$ различное $(n-1)$-буквенное слово.

Введём порядок на этих словах следующим образом: $u\prec v$, если
\begin{itemize}
    \item $u$ лексикографически меньше $v$
    \item представитель, содержащий $u$ левее представителя, содержащего $v$.
\end{itemize}
Из не сильной
$n$-разбиваемости получаем, что максимальное возможное число
попарно несравнимых элементов равно $n-1$. По теореме Дилуорса существует разбиение рассматриваемых $(n-1)$-буквенных слов на $(n-1)$ цепью. Раскрасим слова в $(n-1)$ цвет в соответствии с их принадлежностью к цепям. Раскрасим позиции, с которых начинаются слова, в те же цвета, что и соответствующие слова.

Рассмотрим ориентированный граф $G$ с вершинами вида $(k,i)$, где
$0 < k < n, 0 < i \leqslant l$. Первая координата обозначает цвет, а вторая -- букву.

Ребро с некоторым весом $j$
выходит из $(k_1, i_1)$ в $(k_2, i_2)$, если
\begin{itemize}
    \item для некоторых $i_3,i_4,\ldots,i_{n-1}$ в $j$-ом представителе содержится слово-цикл\linebreak $i_1i_2\cdots i_{n-1},$
    \item позиции, на которых стоят буквы $i_1, i_2$ раскрашены в цвета $k_1, k_2$ соответственно.
\end{itemize}
Таким образом, граф $G$ состоит из ориентированных циклов длины $(n-1)$ с рёбрами одинакового веса. Теперь нам требуется найти показатель, который бы строго монотонно рос с появлением каждого нового представителя при движении от начала к концу слова $W.$ В теореме \ref{verh1} таким показателем было число несократимых рёбер графа $G.$ В доказательстве теоремы \ref{verh2} будет рассматриваться сумма вторых координат неизолированных вершин графа $G.$ Нам потребуется

\begin{lemma}[Основная]\label{n-1:l}
Пусть $A_1, A_2,\ldots,A_{n-1}$ --- вершины графа $G,$ $A_1\to A_2\to\cdots\to A_{n-1}\to A_1$ -- ориентированный цикл длины $(n-1)$ с рёбрами некоторого веса $j.$ Тогда не найдётся другого цикла между вершинами $A_1, A_2,\ldots,A_{n-1}$ одного веса.
\end{lemma}
$\RHD$ От противного. Рассмотрим наименьшее число $j$, для которого нашёлся другой одноцветный цикл между вершинами цикла цвета $j.$ В силу минимальности $j$ можно считать, что этот цикл имеет цвет $k>j.$ Пусть цикл цвета $k$ имеет вид\linebreak $A_{j_1}, A_{j_2},\ldots,A_{j_{n-1}},$ где ${\{j_p\}}_{p=1}^{n-1} = \{1, 2,\ldots,n-1\}.$ Пусть $(k_j, i_j)$ -- координата вершины $A_j.$ Рассмотрим наименьшее число $q\in \mathbb N$ такое, что для некоторого $r$ слово
$i_{j_r}i_{j_{r+1}}\cdots i_{j_{r+q-1}}$ лек\-си\-ко\-гра\-фи\-чес\-ки больше слова $i_{j_r}i_{j_r+1}\cdots i_{j_r+q-1}$ (здесь и далее сложение нижних индексов происходит по модулю $(n-1)$). Такое $q$ существует, так как слова $i_1i_2\cdots i_{n-1}$ и $i_{j_1}i_{j_2}\cdots i_{j_{n-1}}$  сильно сравнимы. Кроме того, в силу совпадения множеств  ${\{j_p\}}_{p=1}^{n-1}$ и $ \{1, 2,\ldots,n-1\}$ получаем, что $q\geqslant 2.$ Так как $q$ --- наименьшее, то для любого $ s<q,$ любого $r$ имеем $ i_{j_r}i_{j_{r+1}}\cdots i_{j_{r+s-1}} = i_{j_r}i_{j_r+1}\cdots i_{j_r+s-1}.$ Тогда для любого $ s<q,$ любого $r$ имеем $ i_{j_{r+s-1}}= i_{j_r+s-1}.$ Из монотонности слов каждого цвета получаем, что для любого $r$ $ i_{j_r}i_{j_{r+1}}\cdots i_{j_{r+q-1}}$ не больше $i_{j_r}i_{j_r+1}\cdots i_{j_r+q-1}.$ Значит для любого $r$ верно неравенство $i_{j_{r+q-1}}\geqslant i_{j_r+q-1}.$ По предположению найдётся такое $r$, что $i_{j_{r+q-1}}> i_{j_r+q-1}.$ Так как обе последовательности ${\{j_{r+q-1}\}}_{к=1}^{n-1}$ и  ${\{j_r+q-1\}}_{к=1}^{n-1}$ пробегают элементы множества чисел $\{1, 2,\ldots,n-1\}$ по одному разу, то $\sum\limits_{r=1}^{n-1}j_{r+q-1} = \sum\limits_{r=1}^{n-1}(j_r+q-1)$ (при вычислении числа $j_r+q-1$ суммирование также проходит по модулю $(n-1)$). Но мы получили $\sum\limits_{r=1}^{n-1}j_{r+q-1} > \sum\limits_{r=1}^{n-1}(j_r+q-1).$ Противоречие. $\LHD$
\medskip

{\bf Завершение доказательства теоремы \ref{verh2}.}

Для произвольного $j$ рассмотрим циклы длины $(n-1)$ цветов $j$ и $j+1$ для некоторого $j.$ Из основной леммы \ref{n-1:l} найдутся числа $k, i$ такие, что вершина $(k,i)$ входит в цикл цвета $(j+1),$ но не входит в цикл цвета $j.$ Пусть цикл цвета $j$ состоит из вершин вида $(k, i_{(j,k)}),$ где $k=1,2,\ldots,n-1.$ Введём величину $\pi(j) = \sum\limits_{k=1}^{n-1}i_{(j,k)}.$ Тогда из основной леммы \ref{n-1:l} и монотонности слов по цветам получаем, что $\pi(j+1)\geqslant\pi(j)+1.$ Так как рассматриваемые периоды не циклические, то найдётся $k$ такое, что $i_{(1,k)}>1.$ Значит, $\pi(1)>n-1.$ $\forall j:i_{(j,k)}\leqslant l-1,$ а значит, $\pi(j)\leqslant (l-1)(n-1).$ Следовательно, $j\leqslant (l-2)(n-1).$ Значит, $t\leqslant (l-2)(n-1).$ Тем самым, теорема \ref{verh2} доказана.

\subsection{Доказательство теоремы \ref{niz}.}

Приведём пример. Из
формулировки этой теоремы следует, что можно положить $l$ сколь
угодно большим. Будем считать, что $l>2^{n-1}$. Мы воспользуемся кон\-струк\-ци\-я\-ми, принятыми в
доказательстве теоремы \ref{verh}. Таким образом, процесс
построения примера сводится к построению рёбер в графе на $l$
вершинах. Разобьём этот процесс на несколько больших шагов. Пусть
на $i$-ом большом шаге в приведённом ниже порядке соединяются
рёбрами следующие пары вершин:\\
$(i,2^{n-2}+i)$,\\
$(i,2^{n-2}+2^{n-3}+i),(2^{n-2}+i,2^{n-2}+2^{n-3}+i)$,\\
$(i,2^{n-2}+2^{n-3}+2^{n-4}+i),(2^{n-2}+i,2^{n-2}+2^{n-3}+2^{n-4}+i),\\
(2^{n-2}+2^{n-3}+i,2^{n-2}+2^{n-3}+2^{n-4}+i),\ldots,$\\
$(i,2^{n-2}+\ldots+2+1+i),\ldots,(2^{n-2}+\ldots+2+i,2^{n-2}+\ldots+2+1+i), $\\
 где $i=2,3,\ldots, l-2^{n-1}+1.$\\
 При этом:
\begin{enumerate}
\item Никакое ребро не будет подсчитано 2 раза, так как вершина
соединена рёбрами только с вершинами, значения в которых
отличаются от значения в выбранной вершине на неповторяющуюся
сумму степеней двойки.

\item Пусть {\it вершина типа} $(k,i)$ --- вершина, которая на $i$-ом
шаге соединяется с $k$ вершинами, значения в которых меньше
значения её самой. Для всех $i$ будут вершины типов $(0,i),
(1,i)\ldots,(n-2,i)$.

Раскрасим в $k$-ый цвет слова, которые для некоторого $i$
начинаются с буквы типа $(k,i)$ и заканчиваются в буквах, с
которыми вершина типа $(k, i)$ со\-е\-ди\-ня\-ет\-ся рёбрами на $i$-ом
большом шаге. Получена корректная раскраска в $(n-1)$ цвет, а,
значит, слово сильно $n$-разбиваемо.

\item На $i$-ом большом шаге осуществляется $(n-2)(n-3)\over 2$
шагов. Значит, $$q=(l-2^{n-1})(n-2)(n-3)/2.$$
\end{enumerate}
Тем самым, теорема \ref{niz} доказана.

\section{Оценка существенной высоты с помощью теоремы~\ref{verh}}
Из рассмотрения случая периодов длины 2 с помощью кодировки букв можно получить оценку на существенную высоту, которая будет расти полиномиально по числу об\-ра\-зу\-ю\-щих и экспоненциально по степени тождества. Для этого надо обобщить некоторые понятия, введённые ранее.

\begin{construction}
Рассмотрим алфавит $\textsuperimposetilde{A}$ с буквами $\{a_1, a_2,\ldots, a_l\}.$ Введём на буквах лексикографический порядок: $a_i>a_j,$ если $i>j.$ Рассмотрим произвольное множество нециклических попарно сильно сравнимых слово-циклов некоторой оди\-на\-ко\-вой длины $t.$ Пронумеруем элементы этого множества натуральными числами, начиная с 1. Введём порядок на словах, входящих в слово-цикл, следующим образом:
$u\prec v,$ если:
\begin{enumerate}
    \item Слово $u$ --- лексикографически меньше слова $v.$
    \item Слово-цикл, содержащий слово $u,$ имеет меньший номер, чем слово-цикл, со\-дер\-жа\-щий слово $v.$
\end{enumerate}
Пронумеруем теперь позиции букв в слово-циклах числами от 1 до $t$ от начала к концу некоторого слова, входящего в слово-цикл.
\begin{notation}
\begin{enumerate}
    \item Пусть $w(i, j)$ -- слово длины $t$, которое начинается с $j$-ой буквы в $i$-ом слово-цикле.
    \item Пусть класс $X(t, l)$ -- рассматриваемое множество слово-циклов с введённым на его словах порядком $\prec.$
\end{enumerate}
\end{notation}

\end{construction}

\begin{definition}
Назовём те классы $X,$ в которых не найдётся антицепи длины $n$, --- {\em $n$-хорошими.} Соответственно те, в которых найдётся такая антицепь --- {\em $n$-плохими.}
\end{definition}

Из теоремы Дилуорса получаем, что слова в $n$-хороших классах $X$ можно рас\-кра\-сить в $(n-1)$ цвет, так что одноцветные слова образуют цепь. Далее требуется оценить число элементов в $n$-хороших классах $X$.

\begin{definition}
Пусть $\beth(t, l, n)$ --- наибольшее возможное число элементов в $n$-хорошем классе $X(t, l).$
\end{definition}
\begin{remark}
Будем считать, что первый аргумент в функции $\beth(\cdot, \cdot, \cdot)$ меньше третьего.
\end{remark}\label{ess:r}
Следующая лемма позволяет оценить $\beth(t, l, n)$ через случаи малых периодов.
\begin{lemma}\label{ess:l1}
$\beth(t, l^2, n)\geqslant\beth(2t, l, n)$
\end{lemma}
$\RHD$ Рассмотрим $n$-хороший класс $X(2t, l).$ Разобьём во всех его слово-циклах позиции на пары соседних так, чтобы каждая позиция попала ровно в одну пару. Затем рассмотрим алфавит $\textsuperimposetilde{B}$ с буквами $\{b_{i,j}\}_{i,j=1}^l,$ причём $b_{i_1, j_1}>b_{i_2,j_2}, $ если\linebreak $i_1\cdot l+j_1>i_2\cdot l+j_2.$ Алфавит $\textsuperimposetilde{B}$ состоит из $l^2$ букв. Каждая пара позиций из разбиения состоит из некоторых букв $a_i, a_j$. Заменим пару букв $a_i, a_j$ буквой $b_{i,j}.$ Поступая так для каждой пары, получаем новый класс $X(t, l^2).$ Он будет $n$-хорошим, так как если в классе $X(t, l^2)$ есть антицепь длины $n$ из слов $w(i_1, j_1), w(i_2, j_2),\ldots,w(i_n, j_n),$ то рассмотрим прообразы слов $w(i_1, j_1), w(i_2, j_2),\ldots,w(i_n, j_n)$ в первоначально взятом классе $X(2t, l).$ Пусть эти прообразы: слова $w(i_1, j'_1), w(i_2, j'_2),\ldots,w(i_n, j'_n)$. Тогда слова $w(i_1, j'_1), w(i_2, j'_2),\ldots,w(i_n, j'_n)$ образуют в классе $X(2t, l)$ антицепь длины $n$. Получено противоречие с $n$-хорошестью класса $X(2t, l).$ Тем самым, лемма доказана.$\LHD$

Теперь оценим $\beth(t, l, n)$ через случаи малых алфавитов.

\begin{lemma}
$\beth(t, l^2, n)\leqslant\beth(2t, l, 2n-1)$
\end{lemma}

$\RHD$ Рассмотрим $(2n-1)$-плохой класс $X(2t, l).$ Можно считать, что $n$ слов из антицепи, а именно, $w(i_1, j_1), w(i_2, j_2),\ldots,w(i_n, j_n)$, начинаются с нечётных позиций слово-циклов. Разобьём во всех его слово-циклах позиции на пары соседних так, чтобы каждая позиция попала ровно в одну пару и первая позиция в каждой паре была нечётной. Затем рассмотрим алфавит $\textsuperimposetilde{B}$ с буквами $\{b_{i,j}\}_{i,j=1}^l,$ причём\linebreak $b_{i_1, j_1}>b_{i_2,j_2}, $ если $i_1\cdot l+j_1>i_2\cdot l+j_2.$ Алфавит $\textsuperimposetilde{B}$ состоит из $l^2$ букв. Каждая пара позиций из разбиения состоит из некоторых букв $a_i, a_j$. Заменим пару букв $a_i, a_j$ буквой $b_{i,j}.$ Поступая так для каждой пары, получаем новый класс $X(t, l^2).$ Пусть слова $w(i_1, j_1), w(i_2, j_2),\ldots,w(i_n, j_n)$ перешли в слова $w(i_1, j'_1), w(i_2, j'_2),\ldots,w(i_n, j'_n).$ Эти слова будут образовывать антицепь длины $n$ в классе $X(t, l^2).$ Таким образом получен $n$-плохой класс $X(t, l^2)$ с тем же числом элементов, что и $(2n-1)$-плохой класс $X(2t, l).$ Тем самым, лемма доказана.$\LHD$

Для дальнейшего рассуждения необходимо связать $\beth(t, l, n)$ для произвольного первого аргумента и для первого аргумента, равного степени двойки.

\begin{lemma}\label{ess:l3}
$\beth(t, l, n)\leqslant\beth(2^s, l+1, 2^s (n-1)+1),$ где $s = \ulcorner\log_2(t)\urcorner.$
\end{lemma}

$\RHD$ Рассмотрим $n$-хороший класс $X(t,l).$ Введём в алфавит $\textsuperimposetilde{A}$ новую букву $a_0$, которая лексикографически меньше любой другой буквы из алфавита $\textsuperimposetilde{A}.$ Получен алфавит $\textsuperimposetilde{A'}.$ В каждый слово-цикл из класса  $X(t,l)$ добавим $(t+1)$-ю, $(t+2)$-ю,$\ldots,$\linebreak $2^s-$ю позиции, на которые поставим буквы $a_0.$ Получили класс $X(2^s, l+1).$ Он будет $(2^s (n-1)+1)$-хорошим, так как в противном случае в этом классе для некоторого $j$ найдутся слова  $w(i_1, j), w(i_2, j),\ldots,w(i_n, j),$ которые образуют антицепь в классе $X(2^s, l+1).$ Тогда
\begin{enumerate}
	\item Если $j>t,$ то слова  $w(i_1, 1), w(i_2, 1),\ldots,w(i_n, 1)$ образуют антицепь в классе  $X(t,l).$
	\item Если $j\leqslant t,$ то слова  $w(i_1, j), w(i_2, j),\ldots,w(i_n, j)$ образуют антицепь в классе  $X(t,l).$
\end{enumerate}
Получено противоречие с тем, что класс $X(t,l)$ --- $n$-хороший. Тем самым, лемма доказана.$\LHD$

\begin{proposition}
$\beth(t, l, n)\leqslant\beth(t, l, n+1)$
\end{proposition}

По лемме \ref{ess:l3} $\beth(t, l, n)\leqslant\beth(2^s, l+1, 2^s (n-1)+1),$ где $s = \ulcorner\log_2(t)\urcorner.$

В силу замечания \ref{ess:r} $t<n$. Значит, $2^s<2n.$

Следовательно, $\beth(2^s, l+1, 2^s (n-1)+1)\leqslant\beth(2^s, l+1, 2n^2).$

По лемме \ref{ess:l1} имеем $$\beth(2^s, l+1, 2n^2)\leqslant\beth(2^{s-1},( l+1)^2, 2n^2)\leqslant\beth(2^{s-2}, (l+1)^{2^2}, 2n^2)\leqslant\beth(2^{s-3}, (l+1)^{2^3}, 2n^2)\leqslant$$$$\leqslant\cdots\leqslant\beth(2, (l+1)^{2^{s-1}},2n^2).$$

По теореме \ref{verh}  имеем $\beth(2, (l+1)^{2^{s-1}},2n^2)<(l+1)^{2^{s-1}}\cdot 4n^4< 4(l+1)^n n^4.$

То есть доказана следующая
\begin{lemma}\label{ess:l4}
$\beth(t, l, n)<4(l+1)^n n^4.$
\end{lemma}

Чтобы применить лемму \ref{ess:l4} к доказательству теоремы \ref{ess:t}, требуется оценить число подслов не $n$-разбиваемого слова с одинаковыми периодами.

\begin{lemma}\label{ess:l5}
Если в некотором слове $W$ найдутся $(2n-1)$ подслов, в которых период повторится больше $n$ раз, и их периоды попарно не сильно сравнимы, то $W$ -- $n$-разбиваемое.
\end{lemma}
$\RHD$ Пусть в некотором слове $W$ найдутся $(2n-1)$ подслов, в которых период повторится больше $n$ раз, и их периоды попарно не сильно сравнимы. Пусть $x$ --- период одного из этих подслов. Тогда в слове $W$ найдутся непересекающиеся подслова $x^{p_1}v'_1,\ldots ,\linebreak x^{p_{2n-1}}v'_{2n-1},$ где $p_1,\ldots ,p_{2n-1}$ --- некоторые натуральные числа, большие $n,$ а $v'_1,\ldots ,v'_{2n-1}$ --- некоторые слова длины $|x|$, сравнимые с $x.$ Тогда среди слов $v'_1,\ldots ,v'_{2n-1}$ найдутся либо $n$ лексикографически больших $x$, либо $n$ лексикографически меньших $x$. Можно считать, что $v'_1,\ldots ,v'_n$ --- лексикографически больше $x$. Тогда в слове $W$ найдутся подслова $v'_1, xv'_2,\ldots ,x^{n-1}v'_n,$ идущие слева направо в порядке лексикографического убывания.$\LHD$

Из этой леммы получаем следствие \ref{co1}.

Рассмотрим не $n$-разбиваемое слово $W.$ Если в нём найдётся подслово, в котором нециклический период $x$ длины $\geqslant n$ повторится больше $2n$ раз, то в слове $x^2$ подслова, которые начинаются с певой, второй$,\ldots,n$-ой позиции, попарно сравнимы. Значит, слово $x^{2n}$ является $n$-разбиваемым. Получаем противоречие с не $n$-разбиваемостью слова $W.$  Из лемм \ref{ess:l5} и \ref{ess:l4} получаем, что существенная высота слова $W$ меньше, чем $(2n-1)\sum\limits_{t=1}^{n-1} \beth(t, l, n) < 8(l+1)^n n^6.$ Значит, $\Upsilon(n, l)<8(l+1)^n n^6.$ Тем самым, теорема \ref{ess:t} доказана.


\begin{thebibliography}{MH1}

\bibitem{Bel04} А.~Я.~Белов. {\it Размерность Гельфанда-Кириллова относительно
свободных ассоциативных алгебр}, Матем. сб., 195:12 (2004), стр.
3--26.

\bibitem{Shirshov1}
Ширшов А.~И. {\it О некоторых неассоциативных ниль-кольцах и
алгебраических алгебрах.} Мат. сб.,
    Т. 41, No 3 (1957), стр. 381--394.

\bibitem{Shirshov2}
 Ширшов А.~И.
 {\it О кольцах с тождественными соотношениями.}
Мат. сб., т. 43, No 2 (1957), стр. 277--283.

\bibitem{Bel07} А.~Я.~Белов. {\it Проблемы бернсайдовского типа, теоремы о высоте и
о независимости}. Фундамент. и прикл. матем., 13:5 (2007), 19--79.


\bibitem{BBL97} А.~Я.~Белов, В.~В.~Борисенко, В.~Н.~Латышев, {\it Мономиальные алгебры.}
Алгебра~--~4, Итоги науки и техн. Сер. Соврем. мат. и ее прил. Темат. обз., 26 (2002), 35--214.

\bibitem{BK} A.J. Belov, M.I. Kharitonov, {\it Subexponential estimations in Shirshov's Height theorem}, arXiv: 1101.4909.

\bibitem{Kuz75} Е.~Н.~Кузьмин. {\it О теореме Нагаты-Хигмана.}
В сб. трудов посвященный 60-летию акад. Илиева. София, 1975. стр.
101--107.

\bibitem{Razmyslov3}
Размыслов Ю.~П. {\it Тождества алгебр и их представлений.}
  {\sl --- М.: Наука, 1989}, 432 стр.


\bibitem{Kem09} Kemer, A.~R., {\it Comments on the Shirshov's Height Theorem.}
    in book: selected papers of A.I.Shirshov, Birkh\"user Verlag AG
    (2009), p. 41--48.

\bibitem{Procesi}
Procesi C. {\it Rings with polynomial identities}, N.Y., 1973, 189
стр.

\bibitem{Ch01} Е.~С.~Чибриков. {\it О высоте Ширшова конечнопорождённой ассоциативной алгебры, удовлетворяющей тождеству степени четыре.} Известия Алтайского государственного университета, 1(19), 2001, стр. 52--56.


\bibitem{Bel92} Belov, A.~Ya., {\it Some estimations for nilpotency of nil-algebras over a field
of an arbitrary characteristic and height theorem}, Commun. Algebra
20 (1992), no. 10, 2919--2922.


\bibitem{BP07} Jean Berstel, Dominique Perrin, {\it The origins of combinatorics on words,} European Journal of Combinatorics 28 (2007) 996–1022.


\bibitem{Lot83} M. Lothaire, {\it Combinatorics of words,} Cambridge mathematical library, 1983.

\bibitem{02} {\it Algebraic Combinatorics on Words}, Cambridge mathematical press, 2002.


\end{thebibliography}
\end{document}